\documentclass[12pt]{article}
\usepackage{epsfig}
\usepackage{subfigure}
\usepackage{authblk}
\usepackage{setspace}
\usepackage{amsmath}
\usepackage{breqn}
\begin{document}
%\FLD{1}{6}{00}{28}{00}

 % Define commands to assure consistent treatment throughout document
 \newcommand{\eqnref}[1]{(\ref{#1})}
 \newcommand{\class}[1]{\texttt{#1}}
 \newcommand{\package}[1]{\texttt{#1}}
 \newcommand{\file}[1]{\texttt{#1}}
 \newcommand{\BibTeX}{\textsc{Bib}\TeX}
 
\title{Numerical anisotropy study of a class of compact schemes}


\author[1]{Adrian Sescu\thanks{{\it email}: sescu@ae.msstate.edu} }  
\author[2]{Ray Hixon }

\affil[1]{Department of Aerospace Engineering, Mississippi State University,  501 Hardy Rd, 330 Walker, Mississippi State, MS 39762}
\affil[2]{MIME Department, University of Toledo, Toledo, OH 43606, USA}

\date{}

\maketitle

\begin{abstract}
We study the numerical anisotropy existent in compact difference schemes as applied to hyperbolic partial differential equations, and propose an approach to reduce this error and to improve the stability restrictions based on a previous analysis applied to explicit schemes. A prefactorization of compact schemes is applied to avoid the inversion of a large matrix when calculating the derivatives at the next time level, and a predictor-corrector time marching scheme is used to update the solution in time. A reduction of the isotropy error is attained for large wave numbers and, most notably, the stability restrictions associated with MacCormack time marching schemes are considerably improved. Compared to conventional compact schemes of similar order of accuracy, the multidimensional schemes employ larger stencils which would presumably demand more processing time, but we show that the new stability restrictions render the multidimensional schemes to be in fact more efficient, while maintaining the same dispersion and dissipation characteristics of the one dimensional schemes.

%\keywords{Partial differential equations \and Wave Propagation \and Compact finite difference schemes \and Stability analysis}
\end{abstract}

\section{Introduction}
\label{intro}

Numerical anisotropy is succinctly defined as the directional dependence of the numerical phase or group velocity, and occurs in numerical approximations of multidimensional hyperbolic partial differential equations. This error is disregarded in many instances on the basis that it does not significantly affect the accuracy of the numerical solution, as the focus is often targeted to other types of errors, such as numerical dissipation, dispersion or aliasing. There are, however, several situations where the numerical anisotropy can deteriorate the numerical solution, as for example in computational acoustics or computational electromagnetics. The isotropy error could be reduced by using, for example, high-resolution spatial schemes or sufficiently dense grids, but the attained reduction may not always suffice. Moreover, by increasing the number of grid points the processing time may increase considerably, while high-resolution spatial schemes may introduce spurious waves at the boundaries.

Most optimizations of finite difference schemes are usually performed in one-dimension, but often the proposed optimized schemes may experience considerably large isotropy errors in multi-dimensions, especially at large wavenumbers. In the past, there were several notable attempts to improve the numerical anisotropy of wave propagation in multi-dimensions. A thorough analysis of the numerical anisotropy was performed by Vichnevetsky~\cite{vich} who solved the two-dimensional wave equation using two different schemes for the Laplacian operator, and averaged the two solutions. A considerable amount of improvement of the isotropy of wave propagation was obtained based on the variation of the weighted average parameters. The same idea was pursued by Trefethen~\cite{trefethen}  who used the leap frog scheme to solve the wave equation in two dimensions. Zingg and Lomax~\cite{zingg1} performed optimizations of finite difference schemes applied to regular triangular grids, that give six neighbor points for a given node. Tam and Webb ~\cite{tam2} proposed an anisotropy correction for Helmholtz equation; they found the anisotropy correction factor applicable to all noise radiation problems, irrespective of the complexity of the noise sources. Lin and Sheu~\cite{lin} used the idea of dispersion-relation-preserving (DRP) of Tam and Webb~\cite{tam1} in two dimensions to optimize the first-order spatial derivative terms of a model equation that resembles the incompressible Navier-Stokes momentum equation. They approximated the derivative using the nine-point grid system, resulting in nine unknown coefficients. Eight of them were determined by employing Taylor series expansions, and the remaining one was determined by requiring that the two-dimensional numerical dispersion relation is the same as the exact dispersion relation. Sun and Trueman~\cite{sun} proposed an optimization of two-dimensional finite difference schemes, by considering additional nodes surrounding the point of differencing. They obtained a significant reduction in the numerical anisotropy, dispersion error and the accumulated
phase errors over a broad bandwidth. Shen and Cangellaris~\cite{shen} introduced a new stencil for the spatial discretization of Maxwell's equations; their scheme experienced superior isotropy of the numerical phase velocity compared with the conventional, second-order accurate FDTD scheme. They also achieved a higher value of the Courant number which implied improved stability, compared to conventional schemes.

In a series of previous papers~\cite{sescu1,sescu2,sescu3,sescu4}, we proposed a technique to derive finite difference schemes in multidimensions with improved isotropy; the proposed explicit schemes employ grid points from more than one direction (as required by the conventional difference schemes). The optimization performed in \cite{sescu1,sescu2,sescu3,sescu4} improved the isotropy of the wave propagation and, moreover, the stability restrictions of the multidimensional schemes in combination with either Runge-Kutta or linear multistep time marching methods. More specifically, it was found that the stability restrictions are more favorable when using multidimensional schemes, even if they involve more grid points in the stencils. However, this was advantageous for low order schemes, such as those of second or fourth order of accuracy, but it was also shown that favorable stability restrictions can be obtained for higher order of accuracy schemes (sixth or eight) by increasing the isotropy corrector factor, while losing the isotropy characteristics. It has been stated in \cite{sescu2} that the time integration was free of errors in order to allow the analysis of the errors resulting from the spatial discretization to be performed. However, that statement was erroneously used since the time marching was not entirely free of errors, but instead the time step was heavily reduced so as the truncation error associated with the time integration can be neglected (we thank one of the reviewer for pointing it).

An attempt to apply the multidimensional optimization to compact schemes is described in this paper, with more focus on stability restrictions. An advantage of using multidimensional compact schemes over explicit schemes (in addition to increased resolution characteristics) is that for a given order of accuracy, fewer stencil points are required. This may be advantageous in the present context, since it will be shown that the optimization applied to compact schemes will provide more convenient stability restrictions, comparable to stability restrictions of low order explicit schemes. It is known that, in general, the computational effort may be higher when using compact schemes since the derivatives are determined in an implicit maner, by inverting large matrices. Hixon~\cite{hixon1} proved that a compact scheme can be split into two backward-forward operators, and the inversion of a tridiagonal matrix can be replaced by the inversion of a bidiagonal matrix. Hixon showed that the calculation of the derivatives can be speeded up since the inversion of a bidiagonal matrix can be done by sweeping from one boundary to the other. Moreover, Hixon~\cite{hixon1} showed that the prefactorization is able to provide higher order of accuracy on a three-point stencil. 

In this paper, it is shown that the speed up can be further increased by using the multidimensional optimization previously applied to explicit schemes~\cite{sescu3,sescu4}. It is found that the speed up for a fourth order accurate prefactored compact scheme can be higher than the speed up provided by a second order explicit multidimensional scheme derived previously, making the multidimensional compact schemes attractive. Moreover, the multidimensional schemes correct the anisotropy of waves propagating in multidimensions (although, for very high order of accuracy, the numerical anisotropy is very low and may be neglected). One may argue, however, that the isotropy error can be corrected by increasing the order of accuracy, but this is not always practical: one may be satisfied with the numerical dispersion and dissipation characteristics of a certain scheme; by employing the multidimensional optimization, the isotropy characteristics and the stability restriction can be improved without affecting the original dispersion and dissipation characteristics. Although is not the subject of this analysis, multidimensional boundary stencils associated with compact schemes can be derived, as was done in \cite{sescu4} for explicit multidimensional schemes.

In section II, the numerical anisotropy is defined briefly, and the multidimensional prefactored compact schemes are introduced and analyzed in spectral space. Details about the numerical anisotropy and the stability restrictions are included. In section III, three test cases consisting of the advection equation in two dimensions, the inviscid Burgers equation in two dimensions, and the advection equations in three dimensions are included and discussed. Conclusions are included in the last section IV.

\section{Multidimensional Prefactored Compact Schemes}
\label{sec:1}

\subsection{Preliminaries}
\label{sec:2}

Consider the initial-value problem in $\textbf{R}^d\times [0,\infty)$:

\begin{equation}\label{ww}
\frac {\partial{u}}{\partial{t}}
=\textbf{c} \nabla u,
\end{equation}

\begin{equation}\label{}
u(\textbf{r},0)=u_0(\textbf{r}),
\end{equation}
where $d$ is the number of dimensions, $\textbf{r}=(x_1,...,x_d)$ is the vector of spatial coordinates, $\textbf{c}=(c_{x_1} \hspace{3mm} ... \hspace{3mm} c_{x_d})$ is the velocity vector, $\nabla=(\partial/\partial{x_1} \hspace{3mm} ... \hspace{3mm} \partial/\partial{x_d})^T$, $u(\textbf{r},t)$ is a scalar function and $u_0(\textbf{r})$ is a periodic function. Let $\Omega=\{(x_1,...,x_d),0<x_1<l_1,...,0<x_d<l_d\}$ be a finite domain in the real space $\textbf{R}^d$ with $l_1$ ... $l_d$ chosen such that there exist a real non-negative number $h=l_1/N_1=...=l_d/N_d$ called space step; $N_1$, ..., $N_d$ are non-negative. For the sake of simplicity, we restrict the analysis to two dimensions ($d=2$); the generalization to higher dimensions is straightforward.

MacCormack~\cite{MacCormack} introduced a two-stage time advancement numerical scheme with a predictor followed by a corrector stage. For the two dimensional advection equation, the schemes can be written as

\begin{eqnarray}\label{e1}
u_{i,j}'   &=& u_{i,j}^{n} - \sigma_x \Delta_{x}^{F} u_{i,j}^{n}
                           - \sigma_y \Delta_{y}^{F} u_{i,j}^{n}  \\
u_{i,j}''  &=& u_{i,j}'    - \sigma_x \Delta_{x}^{D} u_{i,j}'
                           - \sigma_y \Delta_{y}^{D} u_{i,j}'  \\
u_{i,j}^{n+1} &=& \frac{1}{2} (u_{i,j}^{n} + u_{i,j}''),
\end{eqnarray}
where $\sigma_x = c_x k/h $, $k$ is the time step, $\Delta_{x}^{F}$ and $\Delta_{y}^{F}$ are 'forward' difference operators, and $\Delta_{x}^{B}$ and $\Delta_{y}^{B}$ are 'backward' difference operators.

The original scheme is second-order accurate in both space and time and is very well-known due to its simplicity and robustness (Hirsch~\cite{hirsch}). There were attempts to improve the accuracy of the scheme, as was done by Gottlieb and Turkel~\cite{Gottlieb} who increased the spatial accuracy to fourth order, or by Bayliss et al.~\cite{Bayliss} who extended the Gottlieb-Turkel scheme to sixth-order of accuracy.

\subsection{Description and Spectral Analysis of Compact Schemes}

Hixon~\cite{hixon1,hixon2} introduced a new class of high-accuracy compact MacCormack-type schemes for aeroacoustic calculations. These new schemes use a prefactorization to reduce a compact centered difference stencil to two lower-order biased stencils with the weights determined by inverting simpler matrices. The schemes are defined on a three-point stencil and can return up to eight order of accuracy. The general forms of the 'forward' and 'backward' difference operators (along $x$-direction with index $i$) are given by

\begin{eqnarray}\label{e1}
a u_{i+1,j}^{F'} + c u_{i-1,j}^{F'} + (1-a-c) u_{i,j}^{F'} =
\frac{1}{h}
\left[
b u_{i+1,j} - (2b-1) u_{i,j} - (1-b) u_{i-1,j}) 
\right] \nonumber \\
c u_{i+1,j}^{B'} + a u_{i-1,j}^{B'} + (1-a-c) u_{i,j}^{B'} =
\frac{1}{h}
\left[
(1-b) u_{i+1,j} - (2b-1) u_{i,j} - b u_{i-1,j}) 
\right],
\end{eqnarray}
respectively ($j$ is the index along $y$-direction). For fourth order of accuracy,

\begin{eqnarray}\label{e1}
a=\frac{1}{2}-\frac{1}{2\sqrt{3}}; \hspace{6mm} b=1; \hspace{6mm} c=0,
\end{eqnarray}
where the original compact scheme was reduced from one tridiagonal matrix to two bidiagonal matrices which can be solved in an explicit manner. For sixth order accuracy,

\begin{eqnarray}\label{e1}
a=\frac{1}{2}-\frac{1}{2\sqrt{5}}; \hspace{6mm} b=1-\frac{1}{30 a}; \hspace{6mm} c=0,
\end{eqnarray}
where the stencil is reduced from five points to three points with one tridiagonal matrix replaced by two bidiagonal matrices. For eight order of accuracy,

\begin{eqnarray}\label{e1}
a=\frac{1}{4}-\frac{1}{4}\sqrt{\frac{3}{35}} 
-\frac{1}{2}\sqrt{\frac{3}{14}-\frac{1}{2}\sqrt{\frac{3}{35}}};
\hspace{6mm}
b=\frac{70}{69} \left( a^2-\frac{5a}{42} \right);
\hspace{6mm}
c=\frac{1}{70 a},
\end{eqnarray}
and the stencil is reduced from five points to three points with one pentadiagonal matrix replaced by two tridiagonal matrices.

In this work, the fourth and the sixth order accurate prefactored schemes are analyzed since the derivatives can be determined explicitly, by sweeping from one boundary to the other. More suitable expressions of the prefactored schemes (adaptable to the multidimensional optimization) are:

\begin{eqnarray}\label{e1}
u_{i,j}^{F'} = \alpha u_{i+1,j}^{F'} + 
\frac{1}{h}
\left[
b u_{i+1,j} - e u_{i,j}
\right] \nonumber \\
u_{i,j}^{B'} = \alpha u_{i-1,j}^{B'} + 
\frac{1}{h}
\left[
e u_{i,j} - b u_{i-1,j}
\right]
\end{eqnarray}
for fourth order of accuracy, and

\begin{eqnarray}\label{e1}
u_{i,j}^{F'} = \alpha u_{i+1,j}^{F'} + 
\frac{1}{h}
\left[
b u_{i+1,j} - e u_{i,j} - f u_{i-1,j}
\right] \nonumber \\
u_{i,j}^{B'} = \alpha u_{i-1,j}^{B'} + 
\frac{1}{h}
\left[
f u_{i+1,j} - e u_{i,j} - b u_{i-1,j})
\right]
\end{eqnarray}
for sixth order of accuracy, where $\alpha=(a+c-1)/a$, $e=2b-1$ and $f=1-b$ (for fourth order, $f=0$).

\begin{figure}[htp]
\begin{center}
      \label{f1a}\includegraphics[width=0.47\textwidth]{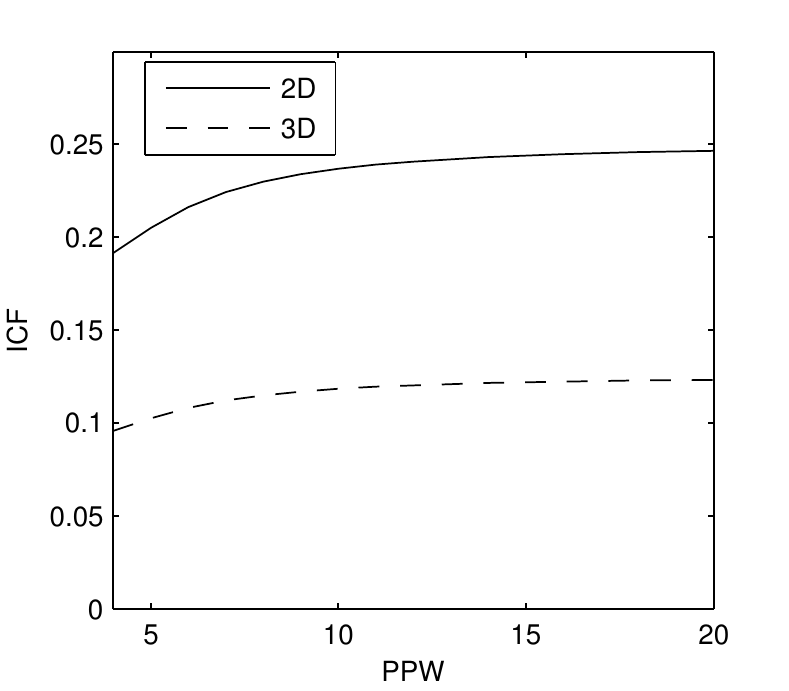}
      \label{f1b}\includegraphics[width=0.47\textwidth]{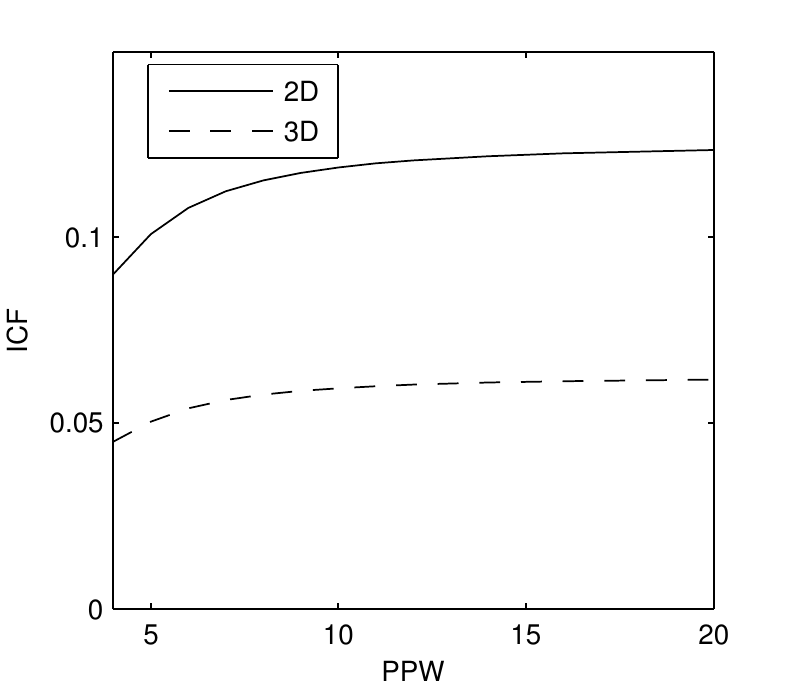}
 \hspace{26mm}(a) \hspace{53mm} (b)
 \end{center}
  \caption{Isotropy corrector factor (ICF) as a function of the number of points per wavelength (PPW): a) fourth order of accuracy; b) sixth order of accuracy.}
  \label{f1}
\end{figure}

Following the same analysis as in~\cite{sescu2}, the multidimensional prefactored compact schemes can be written as

\begin{eqnarray}\label{e1}
u_{i,j}^{F'} &=& \frac{\alpha}{1+\beta} \left[ u_{i+1,j}^{F'} + \frac{\beta}{2} \left( u_{i+1,j-1}^{F'} + u_{i+1,j+1}^{F'} \right) \right]  \\
&+& \frac{1}{h(1+\beta)}
\left[
b u_{i+1,j} - e u_{i,j} 
+ \frac{\beta}{2}
\left(
b u_{i+1,j+1} + b u_{i+1,j-1} - 2e u_{i,j} 
\right)
\right] \nonumber \\
u_{i,j}^{B'} &=& \frac{\alpha}{1+\beta} \left[ u_{i-1,j}^{B'} + \frac{\beta}{2} \left( u_{i-1,j-1}^{B'} + u_{i-1,j+1}^{B'} \right) \right]  \\
&+& \frac{1}{h(1+\beta)}
\left[
e u_{i,j} - b u_{i-1,j}
+ \frac{\beta}{2}
\left(
2e u_{i,j} - b u_{i-1,j+1} - b u_{i-1,j-1}
\right)
\right] \nonumber
\end{eqnarray}
for fourth order of accuracy, and

\begin{eqnarray}\label{e1}
u_{i,j}^{F'} &=& \frac{\alpha}{1+\beta} \left[ u_{i+1,j}^{F'} + \frac{\beta}{2} \left( u_{i+1,j-1}^{F'} + u_{i+1,j+1}^{F'} \right) \right]  \\
&+& \frac{1}{h(1+\beta)}
\left[
b u_{i+1,j} - e u_{i,j} - f u_{i-1,j} 
+ \frac{\beta}{2}
\left(
b u_{i+1,j+1} - f u_{i-1,j-1} +
b u_{i+1,j-1} - f u_{i-1,j+1} - 2e u_{i,j}
\right)
\right] \nonumber \\
u_{i,j}^{B'} &=& \frac{\alpha}{1+\beta} \left[ u_{i-1,j}^{B'} + \frac{\beta}{2} \left( u_{i-1,j-1}^{B'} + u_{i-1,j+1}^{B'} \right) \right]  \\
&+& \frac{1}{h(1+\beta)}
\left[
b u_{i+1,j} - e u_{i,j} - b u_{i-1,j}
+ \frac{\beta}{2}
\left(
f u_{i+1,j+1} - b u_{i-1,j-1} +
f u_{i+1,j-1} - b u_{i-1,j+1} - 2e u_{i,j}
\right)
\right] \nonumber
\end{eqnarray}
for sixth order if accuracy. $\beta$ is the isotropy corrector factor (ICF) and its magnitude, given in figure~\ref{f1} as a function of the number of points per wavelength (PPW), can be determined by minimizing the dispersion error corresponding to the wave-front propagating along a grid line and the wave-front propagating along a diagonal direction (see Sescu et al.~\cite{sescu2} for details). One can notice the the value of ICF in three dimensions is smaller than the value in two dimensions, for the same number of points per wavelength. Apparently, because these multidimensional schemes use a larger number of grid points compared to corresponding one-dimensional schemes, we would expect a larger computational time. Paradoxically, this is not the case because the multidimensional schemes provide more convenient stability restrictions (see Sescu et al.~\cite{sescu3}) as it will be shown in the next section.

Using Fourier analysis, the numerical wavenumbers and the numerical dispersion relation corresponding to the two dimensional wave equation can be found as in~\cite{sescu2}. The numerical dispersion relation has been obtained from the semi-discrete approximation of the wave equation, keeping the time derivative in the continuos form, and discretizing the spatial operator. The individual (forward or backward) numerical wavenumber has both real and imaginary parts: the real part of the forward operator is equal to the real part of the backward operator, and the imaginary parts are opposite. As a result, in a MacCormack predictor-corrector scheme the overall imaginary part is zero. Consequently, the interest is on the real part of the numerical wavenumber. Corresponding to one dimensional prefactored schemes, the real parts of the numerical wavenumbers, $(kh)^*_{4}$ and $(kh)^*_{6}$, are given by the functions

\begin{eqnarray}\label{e1}
Re[(kh)^*_{4}] = f_4(\eta_x) =
\frac
{3 \sin{\eta_x}}
{2+\cos{\eta_x}}
\end{eqnarray}
for fourth order of accuracy, and

\begin{eqnarray}\label{e1}
Re[(kh)^*_{6}] = f_6(\eta_x) =
\frac
{28 \sin{\eta_x} + \sin{2\eta_x}}
{18+12\cos{\eta_x}}
\end{eqnarray}
for sixth order of accuracy, where $\eta_x = k_x h$ and $k_x$ is the $x$-component of the wavenumber. The real parts of the numerical wavenumbers corresponding to multidimensional schemes, for derivatives along $x$-direction, are given by:

\begin{eqnarray}\label{e1}
Re[(kh)^*_{m}] =
\frac{1}{1+\beta}
\left\{
f_m(\eta_x) +
\frac{\beta}{2}
\left[
f_m(\eta_x + \eta_y) +
f_m(\eta_x - \eta_y)
\right]
\right\},
\end{eqnarray}
where $m=4$ for fourth and $m=6$ for sixth order of accuracy, $\eta_y = k_y h$ and $k_y$ is the $y$-component of the wavenumber.

\begin{figure}[htp]
\begin{center}
\label{f2a}\includegraphics[width=0.34\textwidth]{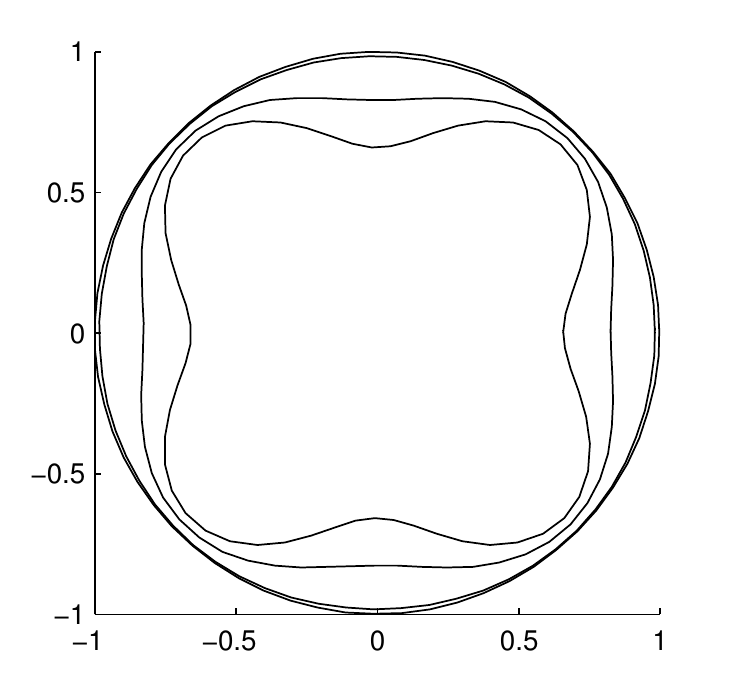}
 \hspace{12mm}
\label{f2b}\includegraphics[width=0.34\textwidth]{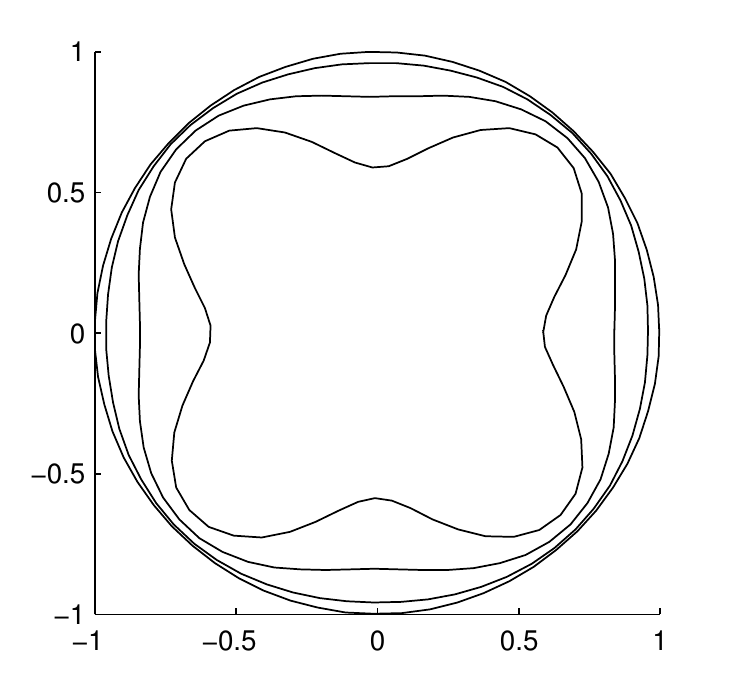}
 \hspace{25mm}(a) \hspace{50mm} (b)
\end {center}
  \caption{Polar diagrams of numerical phase (a) and group (b) velocities as functions of the number points per wavelength, corresponding to one dimensional fourth order schemes.}
  \label{f2}
\end{figure}

\begin{figure}[htp]
\begin{center}
\label{f3a}\includegraphics[width=0.34\textwidth]{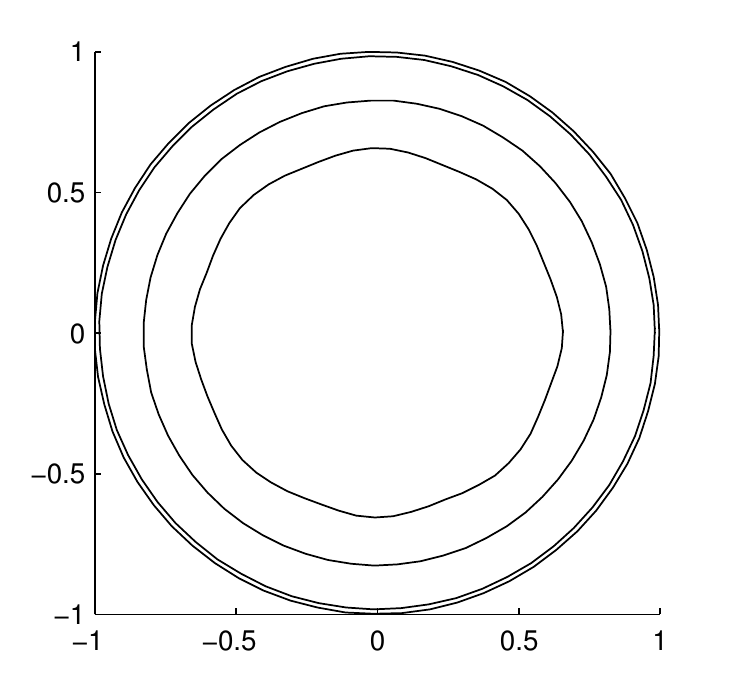}
 \hspace{12mm}
 \label{f3b}\includegraphics[width=0.34\textwidth]{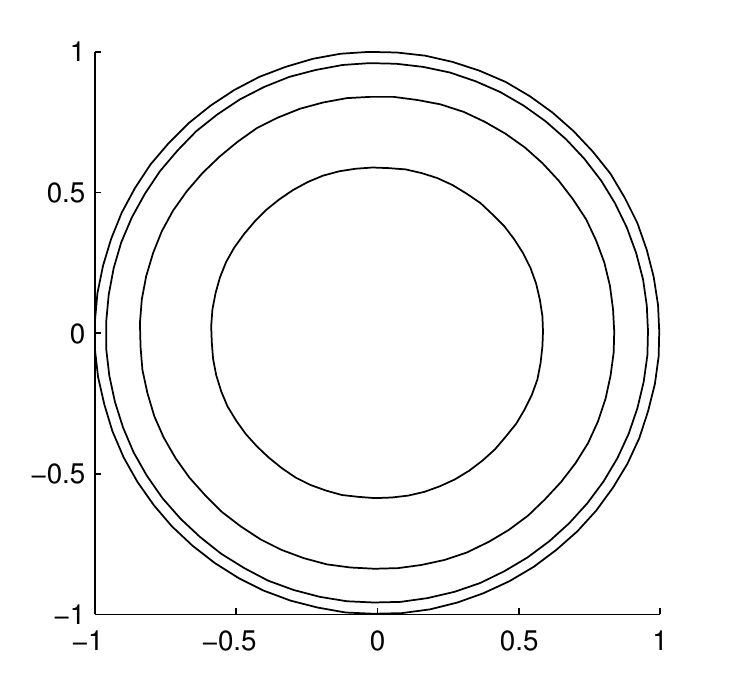}
 \hspace{25mm}(a) \hspace{50mm} (b)
\end{center}
  \caption{Polar diagrams of numerical phase (a) and group (b) velocities as functions of the number points per wavelength, corresponding to multidimensional fourth order schemes.}
  \label{f3}
\end{figure}

\begin{figure}[htp]
\begin{center}
\label{f4a}\includegraphics[width=0.34\textwidth]{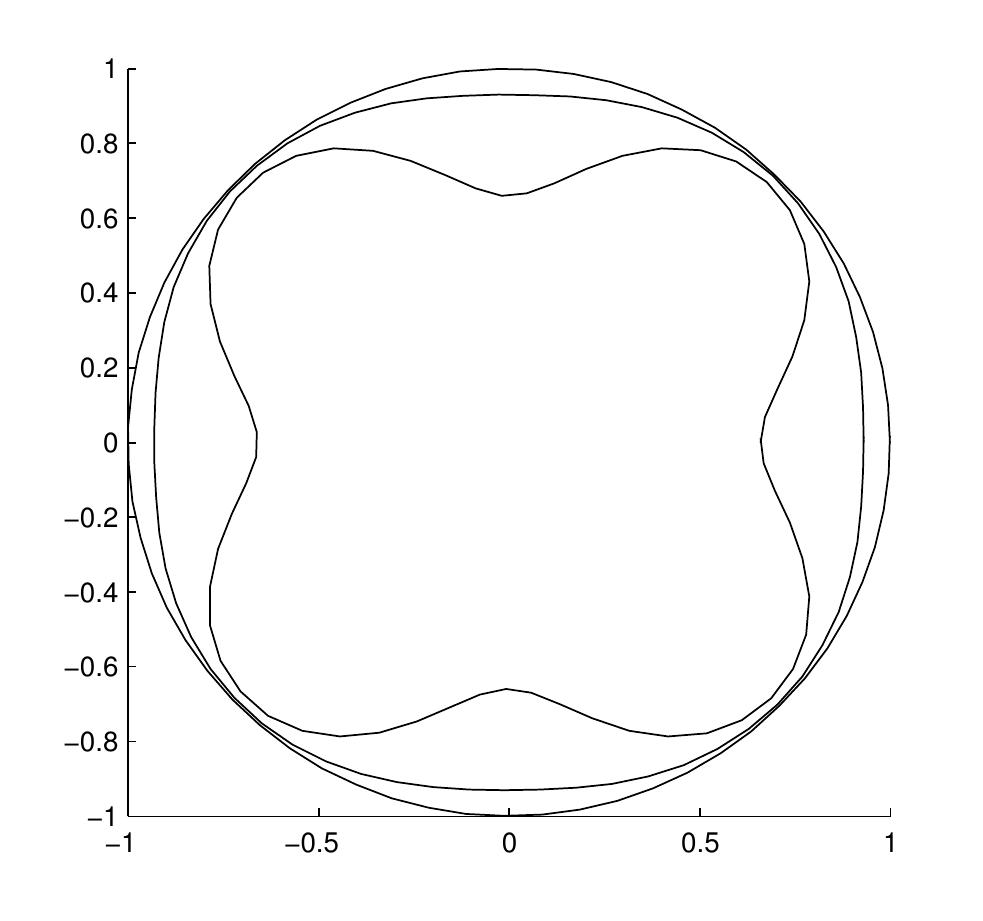}
 \hspace{12mm}
\label{f4b}\includegraphics[width=0.34\textwidth]{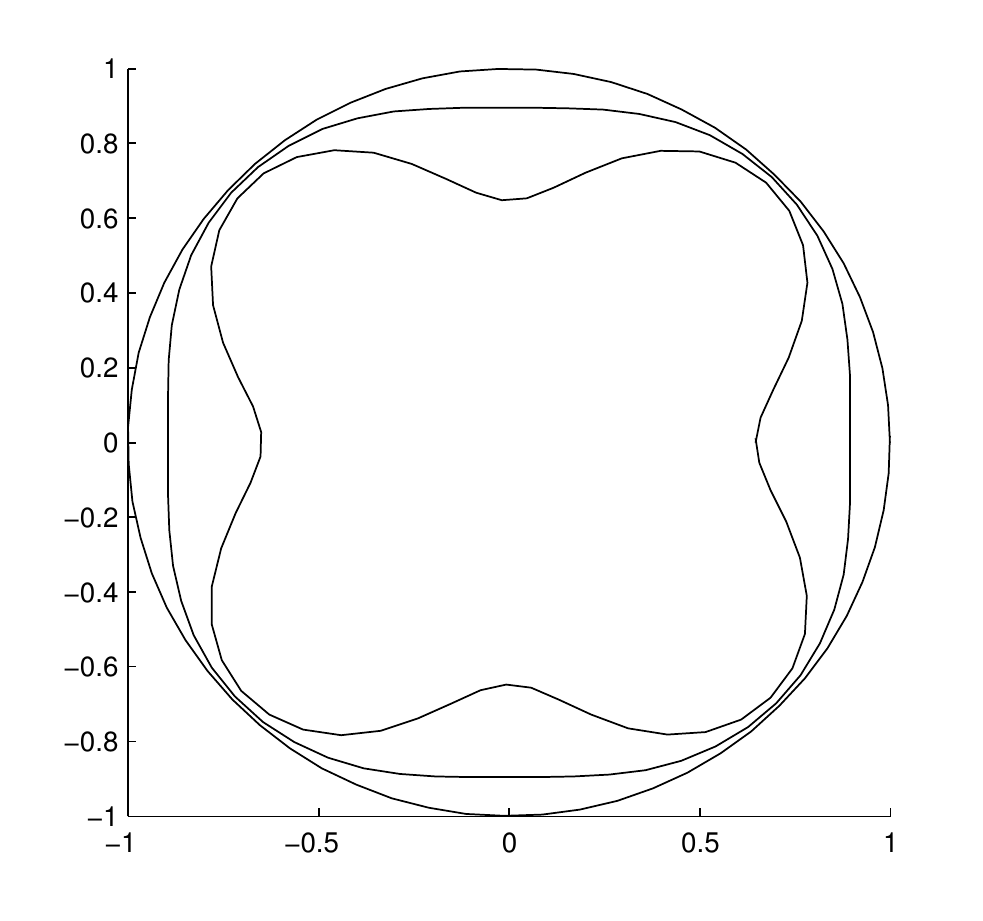}
 \hspace{25mm}(a) \hspace{50mm} (b)
\end{center}
  \caption{Polar diagrams of numerical phase (a) and group (b) velocities as functions of the number points per wavelength, corresponding to one dimensional sixth order schemes.}
  \label{f4}
\end{figure}

\begin{figure}[htp]
\begin{center}
\label{f5a}\includegraphics[width=0.34\textwidth]{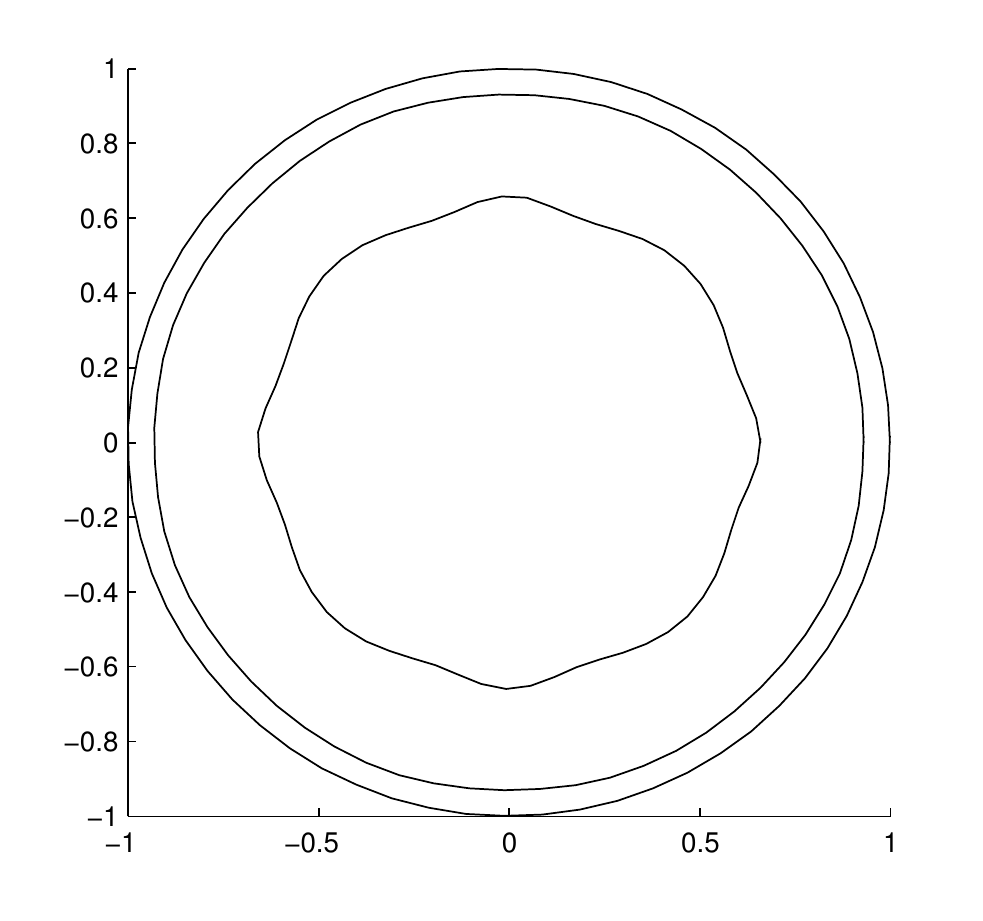}
 \hspace{12mm}
\label{f5b}\includegraphics[width=0.34\textwidth]{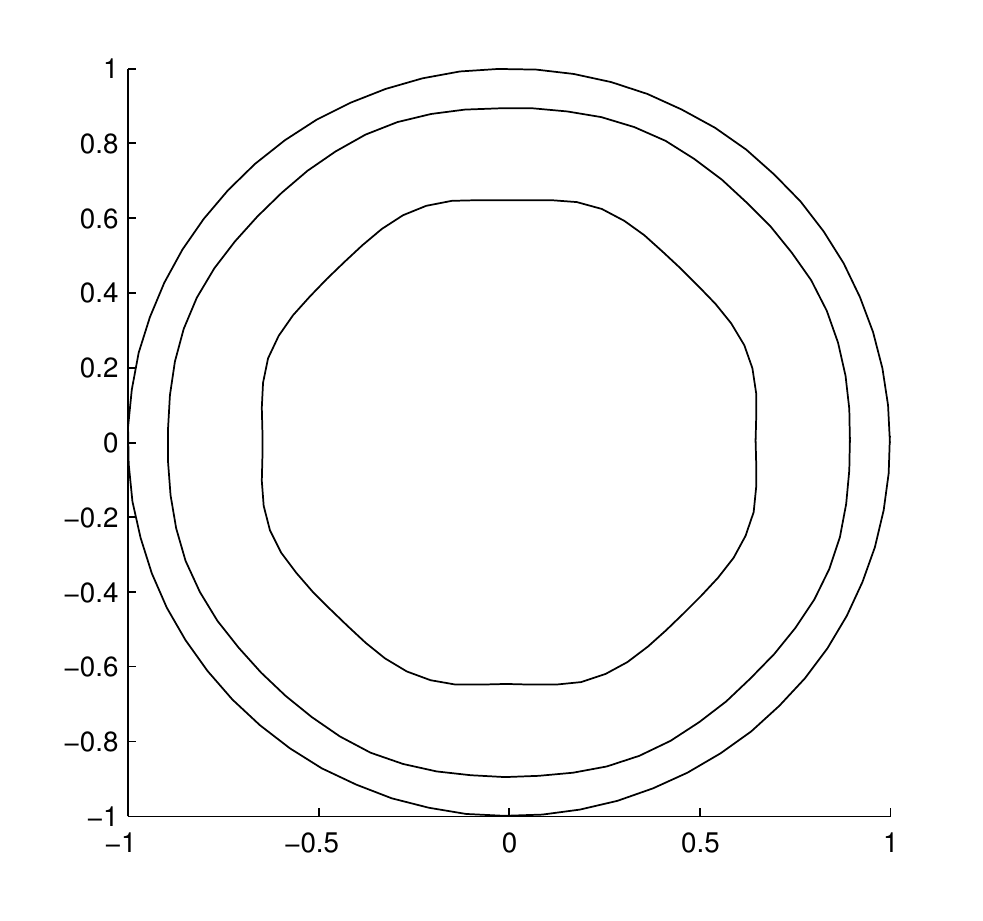}
 \hspace{25mm}(a) \hspace{50mm} (b)
\end{center}
  \caption{Polar diagrams of numerical phase (a) and group (b) velocities as functions of the number points per wavelength, corresponding to multidimensional sixth order schemes.}
  \label{f5}
\end{figure}

Based on the multidimensional numerical dispersion relation, the phase and the group velocities can be calculated. Figure~\ref{f2} show polar diagrams  of numerical phase and group velocities as functions of the number points per wavelength, corresponding to one dimensional fourth order prefactored compact schemes; they show the anisotropy of the schemes. Figure~\ref{f3} show polar diagrams  of numerical phase and group velocities as functions of the number points per wavelength, corresponding to multidimensional fourth order prefactored compact schemes. Figures~\ref{f4} and \ref{f5} show similar plots for sixth order accurate prefactored schemes. Corrections of the numerical anisotropy can be noticed in both figures \ref{f3} and \ref{f5}.

Stability analysis of the multidimensional explicit schemes in combination with Runge-Kutta or linear multistep methods was performed in~\cite{sescu3,sescu4}. The analysis is extended to multidimensional prefactored schemes in the next section, and it will be shown that the multidimensional schemes are more efficient when compared to corresponding one-dimensional schemes.

\subsection{Stability Restrictions}

In this section, stability restrictions imposed by the multidimensional schemes are determined and compared to those imposed by the corresponding one-dimensional schemes, for discretizations of the advection equation in multidimensions.

Let $\Psi_x$ and $-\Psi^{*}_x$ be the Fourier images of the forward and backward differential operators corresponding to $x$-direction, $\Delta_{x}^{F}$ and $\Delta_{x}^{B}$, respectively, and $\Psi_y$ and $-\Psi^{*}_y$ be the Fourier images of the forward and backward differential operators corresponding to $y$-direction, $\Delta_{y}^{F}$ and $\Delta_{y}^{B}$ (the star in $\Psi^{*}_x$ stands for 'complex conjugate'). Also, we denote the Fourier transform of $u^n$ by $\hat{u}^n$, and the amplification factor by $G = \hat{u}^{n+1} / \hat{u}^n$. The Fourier images of the MacCormack stages can be written as

\begin{eqnarray}
\frac{\hat{u}'}{\hat{u}^n} &=& 1 - \sigma_x \Psi_x   - \sigma_y \Psi_y   \\
\frac{\hat{u}''}{\hat{u}'} &=& 1 + \sigma_x \Psi_x^* + \sigma_y \Psi_y^* \\
G &=& \frac{1}{2}
\left(
1 + \frac{\hat{u}''\hat{u}'}{\hat{u}'\hat{u}^n}
\right)
\end{eqnarray}
(see Kennedy and Carpenter~\cite{Kennedy} for one-dimension). To determine the stability restriction, $G$ can be analyzed either for interior schemes or for boundary points using the amplification matrix \textbf{G} (the focus here is on interior points). For stability, the absolute value of $G$ must not exceed $1$. Therefore,

\begin{eqnarray}
|G| = 
\left|
\frac{1}{2}
\left[
1 + (1 - \sigma_x \Psi_x   - \sigma_y \Psi_y)
    (1 + \sigma_x \Psi_x^* + \sigma_y \Psi_y^*)
\right]
\right| \leq 1,
\end{eqnarray}
which gives the stability region for MacCormack scheme. In one dimension the stability restriction for MacCormack scheme is given by $c k/h \leq 1/\xi_{max}$, where $c$ is the advection velocity in one dimension, and $\xi_{max}$ is the maximum numerical wavenumber associated with the spatial differencing scheme used in the predictor and corrector steps (for classical $2$-$2$ MacCormack scheme, $\xi_{max}=1$). In two dimensions, Wendroff~\cite{Wendroff} derived a sufficient condition for the stability region of the classical MacCormack sheme (with $\xi_{max}=1$)) given by

\begin{eqnarray}
\sigma_x^2\sigma_y^2 + \sigma_x^2 + \sigma_y^2 \leq 1/8,
\end{eqnarray}
while Hong~\cite{Hong} derived an exact necessary and sufficient condition given as

\begin{eqnarray}\label{}
\sigma_x^{2/3} + \sigma_y^{2/3} \leq 1.
\end{eqnarray}
The extension of this stability restriction to higher order accurate (in space) MacCormack schemes is given by the equation

\begin{eqnarray}\label{s1}
\sigma_x^{2/3} + \sigma_y^{2/3} \leq \frac{1}{\xi_{max}}.
\end{eqnarray}
where $\xi_{max}$ is approximately $1.73204$ for fourth order, and $1.98943$ for sixth order prefactored compact schemes. For multidimensional MacCormack schemes, the stability resctriction becomes

\begin{eqnarray}\label{r1}
[\sigma_x(1+\beta)]^{2/3} + \sigma_y^{2/3} \leq \frac{(1+\beta)^{2/3}}{\xi_{max}},
\end{eqnarray}
if $|c_x| \geq |c_y|$, and 

\begin{equation}\label{r2}
\sigma_x^{2/3} + [\sigma_y(1+\beta)]^{2/3} \leq \frac{(1+\beta)^{2/3}}{\xi_{max}},
\end{equation}
if $|c_y| \geq |c_x|$. For diagonal directions, with respect to the grid, ($|c_x|=|c_y|=|c|$) the stability restriction becomes

\begin{equation}\label{s2}
\sigma \leq 
\frac{(1+\beta)}{\xi_{max}^{3/2} \left[ 1+(1+\beta)^{2/3} \right]^{3/2}}.
\end{equation}
One can notice that equation (\ref{s1}) along the diagonal direction ($|c_x|=|c_y|=|c|$) becomes $\sigma \leq 1/(2\xi_{max})^{3/2}$, which is consistent with equation (\ref{s2}) when $\beta \rightarrow 0$ (the multidimensional schemes become one-dimensional schemes when $\beta=0$). It is obvious that the right hand side of equation (\ref{s2}) is greater than $1/(2\xi_{max})^{3/2}$ when $\beta > 0$, and goes to $1/(\xi_{max})^{3/2}$ when $\beta \rightarrow \infty$. This generates more efficient stability restrictions for multidimensional schemes. The fact that the prefactored compact shemes require smaller stencils is an advantage when using multidimensional schemes, because the processing time is the result of the balance between the number of additional stencil points and the stability restriction (\ref{s2}).

\section{Test Cases}

\subsection{Advection in a Circular Path}

The rotation of a Gaussian wave in a circular path is considered here to test the stability restrictions derived previously. Consider the advection equation with variable coefficients in $\textbf{R}^2\times [0,\infty)$,

\begin{equation}\label{}
\frac {\partial{u}}{\partial{t}}
=c_x(x,y)\frac {\partial{u}}{\partial{x}}+
 c_y(x,y)\frac {\partial{u}}{\partial{y}}
\end{equation}
with the initial condition,

\begin{equation}\label{o1}
u(x,y,0)=
e^{-ln2
\left(
\frac{(x-x_0)^2+(y-y_0)^2}{\sigma^2}
\right)},
\end{equation}
where 

\begin{equation}\label{o1}
c_x(x,y)=\pi y/2; \hspace{5mm} c_y(x,y)=-\pi x/2
\end{equation}
are the components of the advection velocity, and $\sigma=0.04$. Initially, the Gaussian wave is located in $P(0.25,0)$. Periodic boundary conditions are imposed along both directions.

\begin{table}[htpb]
 \begin{center}
  \caption{\small The number of points in the grids and the corresponding space step for selected spatial schemes ($MC2$ and $MMC2$ - explicit MacCormack schemes of second order; $MC4$ and $MMC4$ - explicit MacCormack schemes of fourth order; $MC6$ and $MMC6$ - explicit MacCormack schemes of sixth order; $PC4$ and $MPC4$ - prefactored compact schemes of fourth order of accuracy; $PC6$ and $MPC4$ - prefactored compact schemes of sixth order of accuracy).}
  \label{t1}
  \begin{tabular}{rrrrrr} \hline
    Spatial scheme &  $N_x \times N_y$ & $h$ \\
\hline
  $MC2$ and $MMC2$ &  $400 \times 400$ & $1.000\times 10^{-2}$ \\
  $MC4$ and $MMC4$ &  $300 \times 300$ & $1.333\times 10^{-2}$ \\
  $MC6$ and $MMC6$ &  $200 \times 200$ & $2.000\times 10^{-2}$ \\
  $PC4$ and $MPC4$ &  $200 \times 200$ & $2.000\times 10^{-2}$ \\
  $PC6$ and $MPC6$ &  $150 \times 150$ & $2.666\times 10^{-2}$ \\
\hline
  \end{tabular}
 \end{center}
\end{table}

The domain is discretized using $N_x \times N_y$ equally-spaced points, where $N_x$ and $N_y$ represent the numbers of grid points along $x$ and $y$ directions, respectively, and are given in table~\ref{t1} for every spatial scheme considered here ($PC4$ and $PC6$ stand for fourth and sixth order accurate conventional prefactored compact schemes, while $MPC4$ and $MPC6$ stand for fourth and sixth order accurate multidimensional prefactored compact schemes, respectively). For comparison in terms of efficiency-gain, numerical results obtained using conventional explicit MacCormack type schemes, of second ($MC2$), fourth ($MC4$), and sixth ($MC6$) order of accuracy, are also included. The corresponding multidimensional explicit schemes are denoted by $MMC2$, $MMC4$, and $MMC6$, respectively. The idea is to compare the gain in efficiency obtained when replacing conventional compact schemes with multidimensional compact schemes with the gain in efficiency when replacing conventional explicit schemes with multidimensional schemes. The time marching is performed using the second order predictor-corrector MacCormack scheme. The values of the ICF are $0.24$ for $MPC4$, and $0.12$ for $MPC6$.

\begin{figure}[htpb]
 \begin{center}
     \label{f6a}\includegraphics[width=0.45\textwidth]{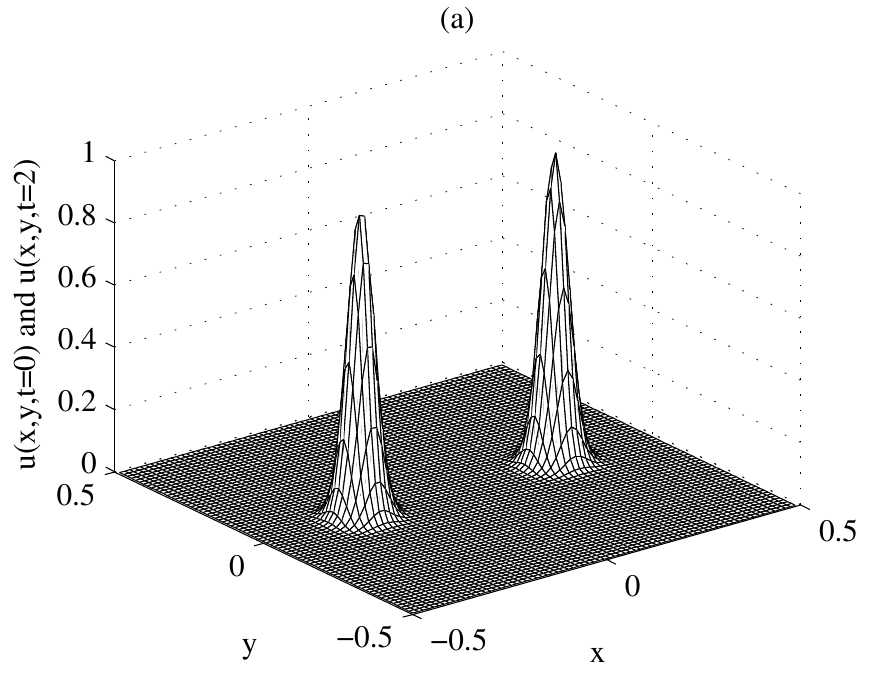}
 \end{center}
 \begin{center}
     \label{f6b}\includegraphics[width=0.45\textwidth]{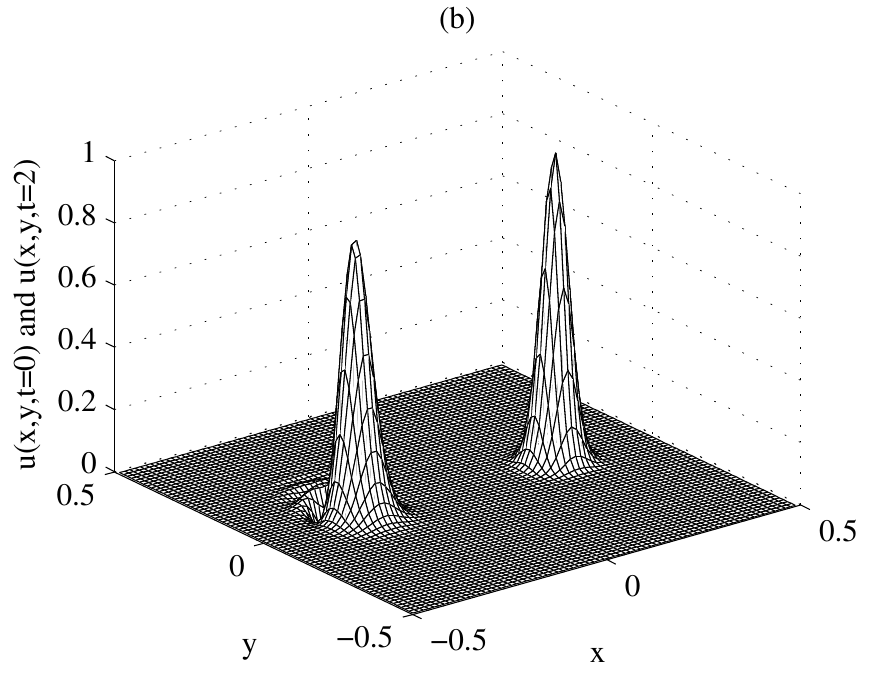}
     \label{f6c}\includegraphics[width=0.45\textwidth]{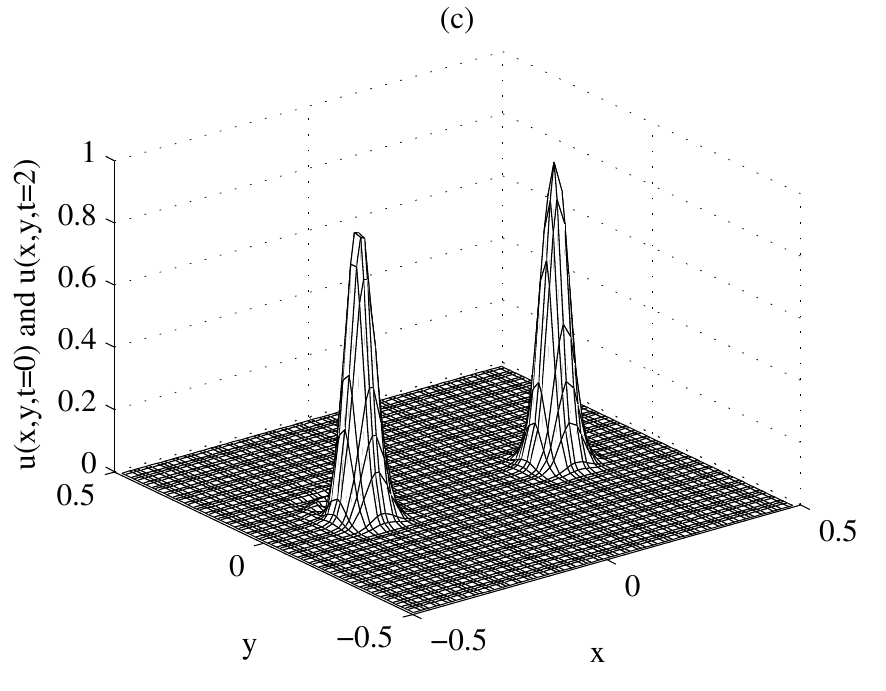}
 \end{center}
  \caption{\small The solution $u(x,y)$ at $t=2$ using multidimensional schemes: a) exact, b) PC4, c) PC6.}
  \label{f6}
\end{figure}

\begin{table}[htpb]
 \begin{center}
  \caption{\small Speed up in percentage by replacing one dimensional schemes with corresponding multidimensional prefactored compact schemes (first column corresponds to isotropy, the second corresponds to maximum allowable time step; $MC2$ and $MMC2$ - explicit MacCormack schemes of second order; $MC4$ and $MMC4$ - explicit MacCormack schemes of fourth order; $MC6$ and $MMC6$ - explicit MacCormack schemes of sixth order; $PC4$ and $MPC4$ - prefactored compact schemes of fourth order of accuracy; $PC6$ and $MPC4$ - prefactored compact schemes of sixth order of accuracy).}
\vspace{2mm}
  \label{t2}
  \begin{tabular}{rrrrrr} \hline
    Stencil &  Isotropy & Max. allowable  \\
\hline
    $MMC2$ compared to $MC2$ &  39\% & 64\% \\
    $MMC4$ compared to $MC4$ &  26\% & 52\% \\
    $MMC6$ compared to $MC6$ &  18\% & 39\% \\
    $MPC4$ compared to $PC4$ &  37\% & 53\% \\
    $MPC6$ compared to $PC6$ &  21\% & 42\% \\
\hline
  \end{tabular}
 \end{center}
\end{table}

The Gaussian (\ref{o1}) revolves around origin for a time interval of $t=2$. Table \ref{t2} shows the speed up in percentage by replacing one dimensional schemes with corresponding multidimensional schemes (for both explicit and compact schemes). Both multidimensional fourth and sixth order accurate are more favorable in terms of the CPU time, compared to corresponding one dimensional schemes. The most important observation when comparing to the results from the explicit schemes, is that the efficiency-gain in terms of the stability restriction obtained when using compact fourth-order schemes is comparable to the efficiency-gain obtained when using explicit second order scheme (and larger when using explicit fourth-order). This is an important results since it shows that by optimizing a compact fourth-order scheme, the same improvement in terms of efficiency is obtained as if second order explicit schemes are optimized. It must be mentioned that the percentage values in table \ref{t2} do not exactly reflect the outcome of equations (\ref{r1}) and (\ref{r2}), since we found by numerical experiments that the actual CFL number can be larger when using multidimensional schemes; the reason lying behind this behavior will be investigated in the future. This is an indication that another more efficient stability restriction must exist, other than equations (\ref{r1}) and (\ref{r2}). The sixth order accurate scheme has a smaller speed up because the isotropy corrector factor is smaller (see figure~\ref{f1}), and in addition more grid points are involved in the stencil. Figure~\ref{f6} shows the initial and the numerical solution $u(x,y)$ at the time $t=2$, multidimensional compact schemes.

\subsection{Inviscid Burgers Equation}

Consider the nonlinear initial-value problem in $\textbf{R}^2\times [0,\infty)$,

\begin{equation}\label{}
\frac {\partial{u}}{\partial{t}}
=u\frac {\partial{u}}{\partial{x}}+
u\frac {\partial{u}}{\partial{y}}
\end{equation}
with the initial condition,

\begin{equation}\label{}
u(x,y,0)=
0.12e^{-
\left(
\frac{x^2+y^2}{\sigma^2}
\right)} + 1.
\end{equation}
The Gaussian is located in $P(0,0)$ and $\sigma=0.2$, with an apmplitude of $0.1$. Periodic boundary conditions are imposed on both directions.

\begin{figure}[htpb]
\begin{center}
\label{f7_a}
     \includegraphics[width=5.8cm]{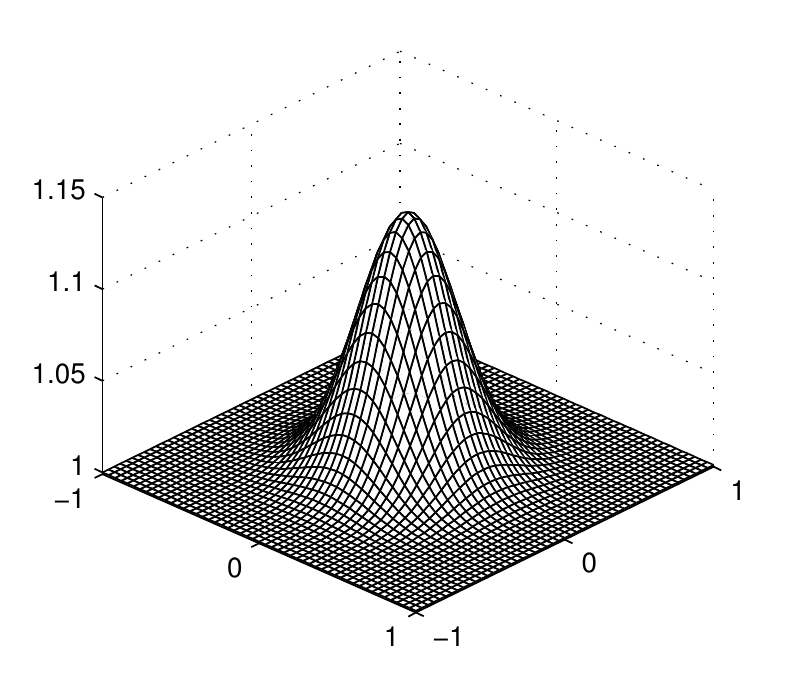}
     \includegraphics[width=5.8cm]{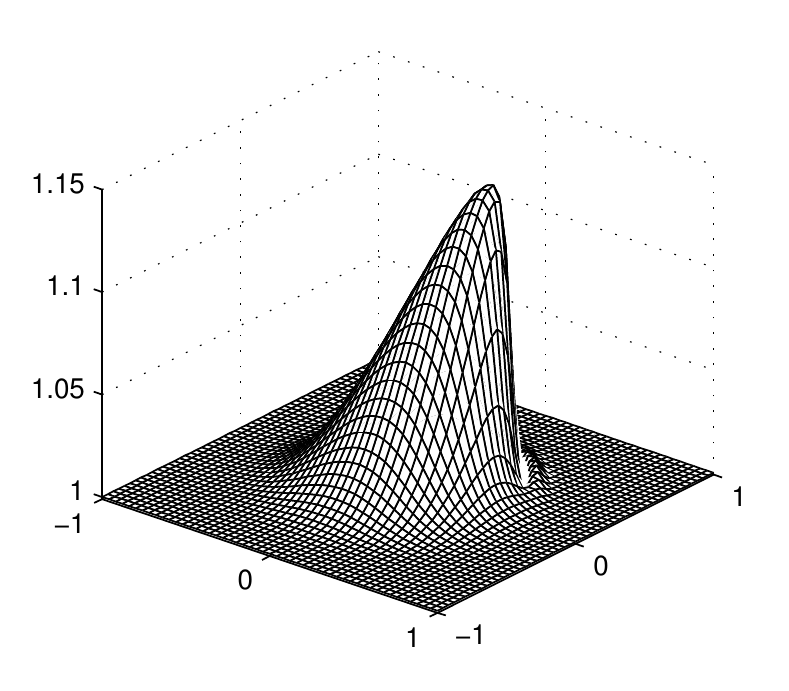}
\end{center}
 \hspace{28mm}(a) \hspace{55mm} (b)
  \caption{\small a) Initial solution for the inviscid Burgers equation b) The solution at time $t=2$s.}
  \label{f7}
\end{figure}

\begin{figure}[htpb]
\begin{center}
\label{f8_a}
     \includegraphics[width=5.8cm]{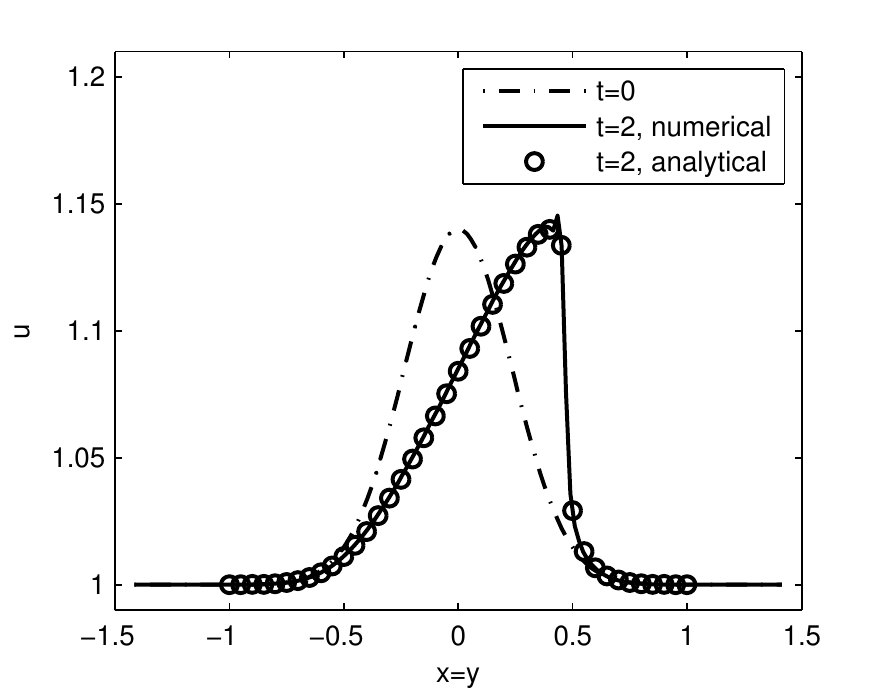}
\end{center}
  \caption{\small Initial solution and the solution to the inviscid Burgers equation at $t=2$s, along the diagonal direction.}
  \label{f8}
\end{figure}

The Gaussian is steepening along the diagonal direction as shown in figure~\ref{f8}. Table~\ref{t3} summarizes the number of points in the grid, corresponding to every spatial scheme(see table \ref{t1} for notations). Figure \ref{f7} shows the initial solution and the numerical solutions obtained using the sixth order pre factored compact schemes, at time $t=2$s, while figure \ref{f8} shows the wave along the diagonal direction, including comparison to the analytical solution. The values of the ICF are $0.24$ for $MPC4$, and $0.12$ for $MPC6$.

\begin{table}[htpb]
 \begin{center}
  \caption{\small The number of points in the grids and the corresponding space step.}
  \label{t3}
  \begin{tabular}{rrrrrr} \hline
    Spatial scheme &  $N_x \times N_y$ & $h$ \\
\hline
  $PC4$ and $MPC4$ &  $150 \times 150$ & $6.666\times 10^{-3}$ \\
  $PC6$ and $MPC6$ &  $100 \times 100$ & $1.000\times 10^{-2}$ \\
\hline
  \end{tabular}
 \end{center}
\end{table}

Table \ref{t4} shows the speed up in percentage by replacing one dimensional schemes with corresponding multidimensional prefactored compact schemes. The same conclusions regarding the efficiency of the multidimensional schemes as for the previous problem apply here. Again the percentages were determined by numerical experiments wherein it was found, as in the previous test case, that the actual CFL number corresponding to multidimensional schemes is larger.

\begin{table}[htpb]
 \begin{center}
  \caption{\small Speed up in percentage by replacing one dimensional schemes with corresponding multidimensional prefactored compact schemes (first column corresponds to isotropy, the second corresponds to maximum allowable time step).}
\vspace{2mm}

  \label{t4}
  \begin{tabular}{rrrrrr} \hline
    Stencil &  Isotropy & Max. allowable  \\
\hline
    $MPC4$ compared to $PC4$ &  26\% & 32\% \\
    $MPC6$ compared to $PC6$&  10\% & 22\% \\
\hline
  \end{tabular}
 \end{center}
\end{table}

\subsection{Advection in Three Dimensions}

Consider the initial-value problem in $\textbf{R}^3\times [0,\infty)$,

\begin{equation}\label{56}
\frac {\partial{u}}{\partial{t}}
=c_x\frac {\partial{u}}{\partial{x}}+
 c_y\frac {\partial{u}}{\partial{y}}+
 c_z\frac {\partial{u}}{\partial{z}}
\end{equation}
with the initial condition,

\begin{equation}\label{57}
u(x,y,z,0)=
e^{-ln2
\left(
\frac{(x-x_0)^2+(y-y_0)^2+(z-z_0)^2}{\sigma^2}
\right)}
\end{equation}
where $c_x$, $c_y$ and $c_z$ are the components of the advection velocity in three dimensions. We choose $c_x=c_y=c_z=1$, $x_0=y_0=z_0=0$ and $\sigma=0.08$, so that the problem can be solved in the domain $[-0.5,0.5]\times [-0.5,0.5]\times [-0.5,0.5]$ with periodic boundary conditions:

\begin{equation}\label{58}
u(0.5,y,z,t)=u(-0.5,y,z,t)
\end{equation}

\begin{equation}\label{59}
u(x,0.5,z,t)=u(x,-0.5,z,t)
\end{equation}

\begin{equation}\label{60}
u(x,y,0.5,t)=u(x,y,-0.5,t)
\end{equation}

Table~\ref{t5} summarizes the number of points in the grid ($N_x \times N_y \times N_z$), corresponding to every spatial stencil. The values of the ICF are $0.11$ for $MPC4$, and $0.055$ for $MPC6$.

\begin{table}[htpb]
 \begin{center}
  \caption{\small The number of points in the grids and the corresponding space step.}
  \label{t5}
  \begin{tabular}{rrrrrr} \hline
    Spatial scheme &  $N_x \times N_y \times N_z$ & $h$ \\
\hline
  $PC4$ and $MPC4$ &  $100 \times 100 \times 100$ & $1.000\times 10^{-2}$ \\
  $PC6$ and $MPC6$ &  $80 \times 80 \times 80$ & $1.25\times 10^{-2}$ \\
\hline
  \end{tabular}
 \end{center}
\end{table}

\begin{table}[htpb]
 \begin{center}
  \caption{\small Speed up in percentage by replacing one dimensional schemes with corresponding multidimensional prefactored compact schemes (first column corresponds to isotropy, the second corresponds to maximum allowable time step).}
\vspace{2mm}

  \label{t6}
  \begin{tabular}{rrrrrr} \hline
    Stencil &  Isotropy & Max. allowable  \\
\hline
    $MPC4$ compared to $PC4$ &  15\% & 25\% \\
    $MPC6$ compared to $PC6$ &  4\% & 14\% \\
\hline
  \end{tabular}
 \end{center}
\end{table}

Table \ref{t6} shows the speed up in percentage by replacing one dimensional schemes with corresponding multidimensional prefactored compact schemes. Both multidimensional fourth and sixth order accurate are more favorable in terms of the CPU time, compared to corresponding one dimensional schemes. However, the percentage values are lower because the isotropy corrector factor is lower in three dimensions, as shown in figure \ref{f1}. The sixth order accurate scheme has a smaller speed up because the isotropy corrector factor is smaller (see figure~\ref{f1}). The CFL corresponding to multidimensional schemes was increased, as before, beyond the limit set by equations (\ref{r1}) and (\ref{r2}).

\section{Conclusions}

This paper has addressed the numerical anisotropy ocurring in finite difference discretizations of hyperbolic partial differential equations, using a predictor-corrector time marching scheme combined with a special class of compact schemes, termed prefactored schemes. It was found that the prefactored compact schemes experience a certain level of numerical anisotropy, despite their high order of accuracy and high resolution characteristics, and that the proposed multidimensional optimization is able to correct the anisotropy. Polar diagrams of numerical phase and group velocities as functions of the number of points per wavelength revealed the improvement in terms of the numerical anisotropy. The isotropy corrector factor has been reported for both two and three dimensional space.

The multidimensional compact schemes were analyzed in terms of numerical stability, and new stability restrictions have been found that make the multidimensional schemes more efficient when compared to corresponding one-dimensional schemes. The stability restrictions given by equations (\ref{r1}) and (\ref{r2}) allow a certain increase of the CFL when multidimensional schemes are employed. However, by numerical experiments we found that the actual CFL can be increased even more when multidimensional schemes are used, resulting in a higher efficiency of the new schemes (as shown in tables \ref{t2}, \ref{t4} and \ref{t4}). The explanation of this behavior will be a subject of future investigations (there must be less restrictive stability limits, other than those given by equations (\ref{r1}) and (\ref{r2})). Three simple test cases showed that the stability restrictions corresponding to the multidimensional compact schemes are more favorable in terms of the processing time.

The multidimensional schemes are suitable to be applied in a parallel code based on domain decomposition. A slight increase in the message passing among the blocks is expected, since not only data along grid lines are exchanged, but also data along diagonal directions. In two dimensions, the additional data to be transferred among the blocks depend on the order of accuracy: for second order explicit or fourth order compact schemes, data associated with one grid point per corner only need to be transferred (four grid points per block); for fourth-order explicit or sixth order compact schemes, data associated with four grid points per corner need to be transferred (sixteen grid points per block). In three dimensions, slighlty more work is involved (data associated with eight corners and 12 sides need to be transferred), but the ratio between the additional data to be transferred when using multidimensional schemes and the data that is transferred when using conventional schemes is comparable to the 2D case. An increase in the processing time per message passing of approximately 2-5\% in 2D and 5-10\% in 3D is expected, depending on the number of grid points of a certain block.

%\section{Acknowledgments}

%The authors would like to express many thanks to the reviewers for their very constructive comments. AS thanks to Dr. A. Afjeh, from the University of Toledo, for many fruitful discussions and for the substantial support.

%\begin{acknowledgements}
%\end{acknowledgements}



\end{document}